\documentclass[11pt,twoside]{article}
\usepackage{latexsym}
\usepackage{amssymb,amsbsy,amsmath,amsfonts,amssymb,amscd,mathtools}
\usepackage{graphicx}
\usepackage{hyperref}
\usepackage{pstricks}
\usepackage[on]{auto-pst-pdf} 
\newcommand{\commentout}[1]{}
\newcommand{\R}{\mathbb{R}}
\def\RR{\mathbb{R}}

\newcommand {\e}  {\varepsilon}
\def\di{\displaystyle}
\newcommand {\da} {\delta}

\newcommand {\Chi} {{\bf \raise 2pt \hbox{$\chi$}} }

\newcommand {\f}   {\frac}
\newcommand {\p}   {\partial}
\newcommand{\fer}{\eqref}
\newcommand{\dis}{\displaystyle}
\newcommand {\proof} {\noindent {\bf Proof}. }
\newcommand{\beq}{\begin{equation}}
\newcommand{\eeq}{\end{equation}}
\newcommand{\bea} {\begin{array}{rl}}
\newcommand{\eea} {\end{array}}
\newcommand{\bepa}{\left\{ \begin{array}{l}}
\newcommand{\eepa} {\end{array}\right.}
\newtheorem{theorem}{Theorem}[section]
\newtheorem{lemma}[theorem]{Lemma}

\newcommand{\qed}{{ \hfill
                       {\unskip\kern 6pt\penalty 500 \raise -2pt\hbox{\vrule\vbox to 6pt{\hrule width 6pt
                       \vfill\hrule}\vrule} \par}   }}
\title{A class of Hamilton-Jacobi equations with constraint: uniqueness and constructive approach
 }

\author{
Sepideh Mirrahimi\thanks{Institut de Math\'ematiques de Toulouse; UMR 5219, Universit\'e de Toulouse; CNRS, UPS IMT, F-31062 Toulouse Cedex 9, France; E-mail: Sepideh.Mirrahimi@math.univ-toulouse.fr} 
\and
Jean-Michel Roquejoffre  \footnotemark[2]\thanks{Institut de Math\'ematiques de Toulouse; UMR 5219, Universit\'e de Toulouse; CNRS, UPS IMT, F-31062 Toulouse Cedex 9, France; E-mail: jean-michel.roquejoffre@math.univ-toulouse.fr}}

\date{\today}

\begin{document}
\maketitle
\pagestyle{plain}
\pagenumbering{arabic}

\begin{abstract}
\noindent We discuss a class of time-dependent Hamilton-Jacobi equations, where an unknown  function of time is intended  to keep the maximum of the solution to the constant value 0. Our main result   is that the full problem has a unique viscosity solution, which is in fact classical.  The motivation is a selection-mutation model which, in the limit of small diffusion, exhibits concentration on the zero level set of the solution of the Hamilton-Jacobi equation. 

\noindent Uniqueness is obtained by noticing that, as a consequence of the dynamic programming principle, the solution of the Hamilton-Jacobi equation is classical. It is then possible to write an ODE for the maximum of the solution, and treat the full problem  as a nonstandard Cauchy problem.
\end{abstract}

\noindent{\bf AMS Class. No:} {35A02, 35F21, 35Q92.}

\bigskip


\section{Introduction}

\subsection{Model and question}
The purpose of this paper is to discuss existence and uniqueness for the following problem, with unknowns $(I(t),u(t,x))$:
\begin{equation}
\label{HJ}
\begin{cases}
\partial_t u = | \nabla  u |^2 +R(x,I)\ (t>0,x\in\R^d), & \di\max_{x} u(t,x) =0\\
I(0)=I_0>0,&\\  
u(0,x)=u_0(x),&
\end{cases}
\end{equation}
where $I_0>0$ and $u_0$ is a concave, quadratic function.
For a given continuous function $I(t)$, $u(t,x)$ solves a  Hamilton-Jacobi equation. The unknown $I$ 
may be thought of as a sort of regulator, or a sort of Lagrange multiplier, to maintain the maximum of $u$ equal to 0.
The constraint on the maximum of $u(t,.)$ makes the problem nonstandard. 

Existence of a solution $(u,I)$  to \eqref{HJ}  is not new, the first result is due to Perthame and Barles \cite{BP} (see also Barles-Perthame-Mirrahimi \cite{BMP} for a result with weaker assumptions). An important improvement is given  
by Lorz, Perthame and the first author in \cite{LMP}; they indeed notice that a concavity assumptions on $R$ - that we also make here - entail regularity. This allows them to
derive the dynamics of the maximum point of a solution $u(t,x)$.  See also \cite{Mthese}. Both types of results rely on a special viscous  approximation of \eqref{HJ} - see equation \eqref{HJE} below.
Uniqueness, however, has remained an open problem, apart from a very particular case \cite{BP}.

The main goal of the paper is to prove the missing uniqueness property; a result that we had already announced in \cite{MR0}. We also provide a constructive existence proof which was not available in the previous existence results \cite{BP,BMP,LMP}. Two important consequences, that we will present in a forthcoming paper \cite{MR2} (see also \cite{MR0}), will be the convergence of the underlying selection-mutation model in a stronger sense than what is known, and asymptotic expansion of the viscous solution.  The asymptotic expansion, which allows to approximate the phenotypical distribution of the population when the mutation steps are small but nonzero, is particularly interesting in view of biological applications. 
 One of the main ingredients will be regularity under suitable concavity assumptions on $R$ and $u_0$, which is far from being available in general. Instead of relying on viscous approximations we will prove these results directly for the  equation 
\begin{equation}
\label{e1.00}
u_t=\vert\nabla u\vert^2+R(t,x).
\end{equation}
This will allow a much easier treatment than in the usual viscosity sense.\\
The uniqueness result will also be helpful to develop the so-called Hamilton-Jacobi approach (see for instance Diekmann et al. 
 \cite{DJMP}, \cite{BP,LMP} and Subsection \ref{sec:mot}) to study more complex models describing selection and mutations. For instance, our result would allow to generalize a result  due to Perthame and the first author in \cite{SM.BP:15} on a selection model with spatial structure, where the proof relies on the uniqueness of the solution to a corresponding Hamilton-Jacobi equation with constraint.

\subsection{Motivation}
\label{sec:mot}

Model (\ref{HJ}) arises in the limit $\varepsilon\to 0$ of the solutions to the problem
\begin{equation} \label{para}
\p_t n_\e - \e \Delta n_\e=  \frac{n_\e}{\e} R \big(x, I_\e(t) \big)\ (t>0, \; x \in \R^d),\  \       \     I_\e(t) = \int_{\R^d}  \psi(x) n_\e(t,x) dx,
\end{equation}
where $n_\e(t,x)$ is the density of a population characterized by a $d$-dimensional biological trait $x$. The population competes for a single resource, this is represented by $I_\e(t)$,
where $\psi$ is a given positive smooth function. The term $R(x,I)$ is the reproduction rate; it is, as can be expected, very negative for large $x$ and decreases as the competition increases. The Laplace term corresponds to the mutations.   The small parameter $\e$ is introduced to consider the long time dynamics of the population when the mutation steps are small. 
Such models can be derived from individual based stochastic
processes in the limit of large populations (see Champagnat-Ferri\`ere-M\'el\'eard \cite{NC.RF.SM:06,CFM}). There is a large literature on the models of population dynamics under  selection and mutations. We refer the interested reader, for instance to Geritz et al. \cite{SG.EK:98} and Diekmann \cite{OD:04} for an approach based on the study of the stability of  differential systems (the so-called adaptive dynamics approach),  to Champagnat \cite{NC.phd} for the study of  stochastic individual based models, to Raoul \cite{Raoulphd} and Mirrahimi \cite{Mthese} for the study of integro-differnetial models.
\\

 The Hopf-Cole transformation $n_\e =\exp \left( {u_\e}/ {\e} \right)$ yields the equation 
\begin{equation}
\label{HJE}
\partial_t u_\e =\e\Delta u_\e +| \nabla  u_\e |^2 +R(x,I_\e)
\end{equation}
 which, in the  limit $\varepsilon\to 0$, yields the equation for $u$. Furthermore, $I_\e$ being uniformly positive and bounded in $\e$, the Hopf-Cole transformation leads to the constraint on $u$.
 One expects that $n_\e$ concentrates at the points where $u$ is close to 0 and the function $I_\e$ appears, in the limit, as a sort of Lagrange multiplier. One has indeed
 $$
n_\e(x,t)\,\underset{\e \rightarrow 0}{\xrightharpoonup{\quad}}\, n(x,t)= \rho(t) \,\delta(x\,-\,\bar x(t)),\quad \text{weakly in the sense of measures},
$$
with
$$
u(t,\overline x(t))=\max_x u(t,x)=0,\qquad \rho(t)= \f{I(t)}{\psi(\overline x(t))}.
$$
 This method to study (\ref{para}) has been introduced in \cite{DJMP} and then developed in different contexts. See for instance Perthame-Barles
  \cite{BP} (convergence  to Hamilton-Jacobi dynamics for  \fer{HJE}, Barles-Mirrahimi-Perthame \cite{BMP} (the same type of result, but with nonlinear, nonlocal diffusion), Champagnat-Jabin \cite{CJ} (nonlinear integro-differential model with several resources), \cite{LMP} (convergence improvement by introduction of the concavity assumptions). This approach has a lot to do with the 'approximation of geometric optics' for reaction-diffusion equations of the Fisher-KPP type. See Freidlin \cite{MFb:85,MF:85} for the probabilistic approach, and  Evans-Souganidis \cite{LE.PS:89}, Barles-Evans-Souganidis \cite{GB.LE.PS:90} for the viscosity solutions approach.


\subsection{Assumptions}
The assumptions  we are stating below are in the same spirit (but slightly weaker) as in \cite{LMP},  where the authors noticed that this set of assumptions allowed them to work with smooth solutions, thus going quite far in the study of \fer{para}. 
We believe that  the results that we will prove   would certainly be false if some of those assumptions  were removed. 
\begin{itemize}
\item{\bf Assumptions on $R(x,I)$.} We choose $R$ to be smooth, and we suppose  that  there is $ I_M>0$  such that (fixing the origin in $x$ appropriately)
\begin{equation} \label{asrmax}
\max_{x \in \R^d} R(x,I_M) = 0 = R(0,I_M) ,
\end{equation}
\begin{equation} \label{asr}
-\underline{K}_1 |x|^2 \leq R(x,I) \leq \overline{K}_0 -\overline{K}_1 |x|^2,  \qquad \text{for }\;  0\leq I \leq I_M,
\end{equation}
\begin{equation} \label{asrD2}
- 2\underline{K}_1 \leq D^2 R(x,I) \leq - 2\overline{K}_1 < 0    \text{ as symmetric matrices},
\end{equation}
\begin{equation} \label{asrDi}
- \underline{K}_2\leq \dis{\frac{\p R}{\p I} \leq - \overline{K}_2},
\end{equation}
\beq
\label{asr23}
| \f{\p^2 R}{\p I \p x_i}(x,I)| + | \f{\p^3 R}{\p I \p x_{i} \p x_j}(x,I)| \leq K_3,  \qquad \text{for  $0\leq I \leq I_M$, and $i,\,j=1,2,\cdots,d$},
\eeq
\begin{equation}\label{asRD3}
\|D^3R(\cdot,I) \|_{L^\infty(\R^d)}\leq K_4, \qquad \text{ for }0\leq I \leq I_M.
\end{equation}
\item{\bf Assumptions on $u_0(.)$ and $I_0$.}
We assume the existence of positive constants $\underline{L}_0, \overline{L}_0,\underline{L}_1, \overline{L}_1$ such that
\begin{equation} \label{asu}
-\underline{L}_0 -\underline{L}_1 |x|^2 \leq u_0(x) \leq  \overline{L}_0 - \overline{L}_1 |x|^2 ,
\end{equation}
\begin{equation} \label{asuD2}
-2\underline{L}_1 \leq D^2u_0  \leq - 2\overline{L}_1 .
\end{equation}
Note that this implies
\begin{equation}
\label{asuD1}
\vert Du_0(x)\vert\leq L_2(1+\vert x\vert),
\end{equation}
for a large constant $L_2>0$.
We also need that, for a positive constant $L_3$,
\begin{equation} \label{asuD3}
\|D^3u_0 \|_{L^\infty(\R^d)} \leq L_3.
\end{equation}
Finally we assume that
\beq
\label{as:u0-I0}
\max_x \, u_0(x) = u_0 (\overline x_0)=0, \qquad R(\overline x_0, I_0)=0.
\eeq

\end{itemize}

Note that the monotony assumption \fer{asrDi} means that the growth rate decreases as the competition increases, which is natural from the modeling point of view.  The concavity assumption is a technical one.\\

In Section \ref{sec:pre} we will study an unconstrained Hamilton-Jacobi equation where we replace $R(x,I)$ by $R(t,x)$. To prove our results on this unconstrained problem we assume same type of regularity and concavity assumptions on $R$ that we state below:
\begin{itemize}
\item{\bf Assumptions on $R(t,x)$.} We choose $R$ to be smooth, and we suppose  that 
\begin{equation} \label{asr-t}
-\underline{K}_1 |x|^2 \leq R(t,x) \leq \overline{K}_0 -\overline{K}_1 |x|^2,  \qquad \text{for }\;  t\in \R^+,
\end{equation}
\begin{equation} \label{asrD2-t}
- 2\underline{K}_1 \leq D^2 R(t,x) \leq - 2\overline{K}_1 < 0    \text{ as symmetric matrices},
\end{equation}
\begin{equation}\label{asRD3-t}
\|D^3R(t,\cdot) \|_{L^\infty(\R^d)}\leq K_4, \qquad \text{ for }t \in \R^+.
\end{equation}
\end{itemize}
\subsection{Results and plan of the paper}

Our first result concerns the unconstrained Hamilton-Jacobi equation
\beq
\label{e1.1}
\begin{cases}
u_t=\vert\nabla u\vert^2+R(t,x) &   (t>0,x\in\R^d),\\
u(0,x)=u_0(x).&
\end{cases}
\eeq
The assumptions on $u_0$ and $R$ are those stated in the preceding subsection. 
\begin{theorem}
\label{cauchy}
Equation \fer{e1.1} has a unique viscosity solution $u$ that is bounded from above. Moreover, it is a classical solution: $u\in L^\infty_{\mathrm loc} \big(\R^+; W^{3,\infty}_{\mathrm loc}(\R^d) \big)\cap W^{1,\infty}_{\mathrm loc} \big(\R^+; L^{\infty}_{\mathrm loc}(\R^d) \big)$, $- \max( 2 \, \underline L_1,\sqrt{\underline K_1}) \leq D^2 u \leq - \min (2\overline L_1,\sqrt{\overline K_1})$ and $\| D^3 u\|_{L^\infty([0,T] \times \R^d)} \leq L_4(T)$ where $L_4(T)$  is a positive constant depending on $\underline L_1,\, \underline K_1,\, K_4, \,L_3$ and $T$. 
\end{theorem}
Let us point out that the assumption of  boundedness from above is, most certainly, irrelevant. However the constraint in \fer{HJ} implies that all the solutions that we consider are bounded from above. This extra condition is thus legitimate, and will keep the length of the preliminary work to  a minimum.
\begin{theorem}
\label{thm:uniq}
The Hamilton-Jacobi equation with constraint \fer{HJ} has a unique solution $(u,I)$. Moreover we have
$$
(u,I)\in L^\infty_{\mathrm loc} \big(R^+; W^{3,\infty}_{\mathrm loc}(\R^d) \big)\cap W^{1,\infty}_{\mathrm loc} \big(R^+; L^{\infty}_{\mathrm loc}(\R^d) \big) \times W^{1,\infty}(\R).
$$
\end{theorem}

\bigskip

The paper is organised as follows. In Section 2 we prove the Cauchy Problem for \eqref{e1.1}. In Section 3 we reduce \eqref{HJ} to a (nonstandard) differential system. Theorem \ref{thm:uniq} is proved in Section 4. Section 5 is devoted to the study of a particular example.

\section{The Cauchy problem}
\label{sec:pre}
In this section, we prove Theorem \ref{cauchy}. \\

\noindent
We will first prove that  the only solution to \eqref{e1.00} that is bounded from above is the solution $u(t,x)$ of the dynamic programming principle
\beq
\label{dyn-p}
u(t,x) = \sup_{\underset{\gamma(t)=x}{(\gamma(s),s)\in \R^d\times[0,t]}} \left\{ F(\gamma) \ :\  \gamma \in C^1([0,t];\R^d) \right\},
\eeq
with
$$
F(\gamma):= \ u_0(\gamma (0)) + \int_0^t \left( -\f{|\dot{\gamma}Ê|^2}{4}(s)+R(s,\gamma(s)) \right) ds.
$$
Such a solution is a viscosity solution to \fer{HJ}. We will prove, in addition, that it is classical and satisfies the properties claimed by Theorem \ref{cauchy}.
\\

\noindent{\bf Uniqueness.} This step essentially consists in showing that a viscosity solution of \fer{e1.00} does not grow too wildly, which will reduce the problem to the application of classical arguments. We have already assumed boundedness from above, so let us show that $u$ goes to $-\infty$ at most in a quadratic fashion. Due to the assumptions  \fer{asr-t} on $R$, we have
$$
\partial_tv\geq 0, \   \  v(t,x)=u(t,x)+t\underline K_1\vert x\vert^2
$$
in the viscosity sense, which implies that $v$ is time-increasing - thus the needed estimate. Let us - although this is elementary - explain why: choose $T>0$ and assume the existence of $0<s<t\leq T$ such that the inequality $v(t,x)\geq v(s,x)$ does not hold. In other words there is $x_0$ such that $v(s,x_0)>v(t,x_0)$. For $\e>0$ consider the quantity
$$
w_\e(t,x):=v(t,x)-v(s,x_0)+\frac\e{T-t}+\frac{\vert x-x_0\vert^2}{\e^2}.
$$
For $\e>0$ small enough, there is a local minimum point for $w_\e$, called $(t_\e,x_\e)$, such that $t_\e$ is bounded away from $s$ and $T$, and $x_\e\to x_0$ as $\e\to0$. At that minimum point the viscosity inequality implies $-\di\partial_t\frac1{T-t}\geq0$, a contradiction.

So there is a large $C>0$ such that 
$$
-C(1+t)(1+\vert x\vert^2)\leq u(t,x)\leq C.
$$
And so, by an easy adaptation of Chap. 2 of  Barles \cite{Barles},  where a uniqueness result for solutions that grow at most exponentially fast is provided (see also \cite{FIL1,FIL2}), there is at most one viscosity solution to \fer{e1.00} that is bounded from above.\\

\noindent {\bf Existence.} We may thus turn to \fer{dyn-p}. Let us suppose that $(\gamma_n)_{1 \leq n}$, with $\gamma_n \in C^1([0,t];\R^d)$ and $\gamma_n(t)=x$,   is such that $F(\gamma_n) \to u(t,x)$ as $n\to \infty$. Since $R$ and $u_0$ are bounded from above, we obtain that, for some constant $C$
$$
\int_0^t |\dot{\gamma}_n |^2(s)ds < C.
$$
Consequently, from  $\gamma_n(t)=x$ we deduce, modifying the constant $C$ if necessary, that 
$$
\| \gamma_n \|_{W^{1,2}[0,t]} <C.
$$
It follows that, there exists $\overline \gamma\in W^{1,2}([0,t] ; \R^d)$, such that as $n \to \infty$, $\gamma_n \to \overline \gamma$ strongly in $C([0,t]; \R^d)$ and weakly in $W^{1,2}([0,t] ; \R^d)$. We deduce that, as $n\to \infty$,
$$
u_0(\gamma_n(0))\to u_0(\overline \gamma), \quad \int_0^t R(s,\gamma_n(s))ds \to \int_0^t R(s,\overline \gamma(s))ds, \quad
 \int_0^t |\dot{\overline \gamma} |^2(s)ds \leq \liminf_{n\to \infty}  \int_0^t |\dot{\gamma}_n |^2(s)ds.
$$
We conclude that
\begin{equation}
\label{e2.2}
u(t,x) = \ u_0(\overline \gamma (0)) + \int_0^t \left( -\f{|\dot{\overline \gamma}|^2}{4}(s)+R(s,\overline \gamma(s)) \right) ds.
\eeq
We claim that such a trajectory is unique, which entails uniqueness for the Cauchy problem \eqref{e1.1}. We note indeed that such trajectory $\overline \gamma$ satisfies the following Euler-Lagrange equation 
\beq
\label{E-L}
\begin{cases}
\ddot{\overline \gamma}(s)=-2 \nabla R(s,\overline \gamma(s)),\\
\dot{\overline \gamma }(0)=-2\nabla u_0(\overline  \gamma (0)),\\
\overline \gamma(t)=x.
\end{cases}
\eeq
However, from the concavity assumptions on $R$ and $u_0$ we obtain that the  above elliptic problem is  coercive  and hence the solution $\overline \gamma$ is unique.\\

\noindent {\bf Regularity}. Let us denote by $\gamma_x(t)$ the unique solution of \eqref{E-L}. The function $(t,x)\mapsto(\gamma_x(t),\dot\gamma_x(t))$ belongs to 
$L^\infty_{\rm loc}(\RR^+,W^{2,\infty}_{\rm loc}(\RR^{d}))$ because $\nabla u_0$  and $\nabla R(t,.)$ are in $W^{2,\infty}_{\rm loc}(\RR^{d})$. Now, we have
$$
u(t,x)=u_0(\gamma_x(0))+F(\gamma_x),
$$
yielding (this is a classical computation):
$$
\partial_iu(t,x)=\nabla u_0(\gamma_x(0)).\partial_i\gamma_x(0)+\int_0^t \left( -\, \frac{\dot\gamma_x(s).\partial_i\dot\gamma_x(s)}{2}+\nabla R(s,\gamma_x(s))\partial_i\gamma_x(s)\right)ds.
$$
Integrating by parts and using the Euler-Lagrange equation \fer{E-L}, we get:
\beq
\label{e2.4}
\nabla u(t,x)=-\frac{\dot\gamma_x(t)}2.
\eeq
Thus $u\in L^\infty_{loc}(\RR^+,W^{3,\infty}_{loc}(\RR^d)).$
\\

\noindent {\bf Strict concavity.} We will prove that $u(t,x)$ is uniformly strictly concave, namely that $D^2u\leq -2 \lambda I$ in the sense of symmetric matrices, for $\lambda = \min \big(\overline L_1,\f{\sqrt{\overline K_1}}{2} \big)$. To this end, we show 
that, for all $\sigma\in[0,1]$ and $(x,y)\in\RR^d\times\RR^d$:
\beq
\label{conc-u}
\sigma u(t,x) + (1-\sigma) u(t,y) +\lambda \sigma (1-\sigma) |x-y|^2 \leq u(t, \sigma x+(1-  \sigma ) y).
\eeq
Let $\gamma_x$ and $\gamma_y$ be optimal trajectories, solving \fer{E-L}, with 
$\gamma_x(t)=x$ and $\gamma_y(t)=y$.
%
Note from the choice of $\gamma_x$ and $\gamma_y$ that we have
$$
u(t,x)= u_0(\gamma_x(0) )+  \int_0^t \left( -\f{|\dot{ \gamma}_x(s)\vert^2}4+R(s, \gamma_x(s)) \right) ds,
$$
$$
u(t,y)= u_0(\gamma_y(0) )+  \int_0^t \left( -\f{|\dot{ \gamma}_y(s)\vert^2}4+R(s, \gamma_y(s)) \right) ds,
$$
and
$$
\begin{array}{rll}
u(t,\sigma x+(1-\sigma)y)\geq &u_0(\sigma\gamma_x(0)+(1-\sigma)\gamma_y(0))\\
&+\di\int_0^t\left(-\f{|\sigma\dot{ \gamma}_x + (1-\sigma)\dot{ \gamma}_y(s)\vert^2}4+R(s,\sigma\gamma_x(s)+(1-\sigma)\gamma_y(s))\right)ds.
\end{array}
$$
Furthermore, from the concavity assumptions on $R$ and $u_0$ we have
$$
\sigma u_0(t,x) + (1-\sigma) u_0(t,y) + \overline L_1 \sigma (1-\sigma) |\gamma_x(0)-\gamma_y(0)|^2 \leq u_0(t, \sigma x+(1-  \sigma ) y),
$$
and
$$
\begin{array}{c}
\sigma  \di\int_0^t  R(s, \gamma_x(s))  ds + (1-\sigma)   \di\int_0^t  R(s, \gamma_y(s))  ds + \overline K_1 \sigma (1-\sigma)\di\int_0^t |\gamma_x(s)-\gamma_y(s)|^2ds \\
\leq \di\int_0^t  R(s, \sigma \gamma_x(s) + (1-\sigma) \gamma_y(s))  ds.
\end{array}
$$
Moreover, from the strict concavity of $\mu\mapsto -|\mu|^2$, we obtain that
$$
\begin{array}{c}
\sigma \di\int_0^t  -\f{|\dot{ \gamma}_x(s)\vert^2}4 ds +(1-\sigma)  \di\int_0^t  -\f{|\dot{ \gamma_y}(s)\vert^2}4 ds+\sigma(1-\sigma) \di\int_0^t  \f{|\dot{ \gamma}_x (s)- \dot{ \gamma}_y(s)\vert^2}4 ds
\\
\leq \di\int_0^t  -\f{|\sigma\dot{ \gamma}_x + (1-\sigma)\dot{ \gamma}_y(s) \vert^2}4ds.
\end{array}
$$
We deduce that
\beq
\label{e2.8}
\begin{array}{rll}
&u(t,\sigma x+(1-\sigma)y))\geq\sigma u(t,x)+(1-\sigma)u(t,y)\\
+&\sigma(1-\sigma)\left(\di\di\int_0^t(\f 1 4  \vert\dot\gamma_x(s)-\dot\gamma_y(s)\vert^2+ \overline K_1 \vert\gamma_x(s)-\gamma_y(s)\vert^2)ds+\overline L_1 \vert\gamma_x(0)-\gamma_y(0)\vert^2\right)
\end{array}
\eeq
Next we have
$$
 \f{\sqrt{\overline K_1}}{2}    \int_0^t \f{d}{ds} |\gamma_x(s)-\gamma_y(s)|^2ds \leq   \overline K_1 \int_0^t |\gamma_x(s)-\gamma_y(s)|^2ds +\int_0^t  \f{|\dot{ \gamma}_x - \dot{ \gamma}_y|^2}{4}(s) ds.
$$
Writing
$$
\vert x-y\vert^2=\vert\gamma_x(t)-\gamma_y(t)\vert^2=\vert\gamma_x(0)-\gamma_y(0)\vert^2+\int_0^t\frac{d}{ds}\vert\gamma_x(t)-\gamma_y(t)\vert^2 ds
$$
we find
$$
\, \vert x-y\vert^2\leq    \vert\gamma_x(0)-\gamma_y(0)\vert^2+    2\sqrt{\overline K_1} \int_0^t |\gamma_x(s)-\gamma_y(s)|^2ds +
\f{1}{2\sqrt{\overline K_1}}\int_0^t  |\dot{ \gamma}_x - \dot{ \gamma}_y|^2(s) ds.
$$
Combing the above line with \fer{e2.8}, we obtain    \fer{conc-u} for $\lambda=\min \big(\overline L_1,\di\f{\sqrt{\overline K_1}}{2}\big)$.
\\

\noindent{\bf Bounds on $u$ and $\nabla u$.}  The first thing to notice is a bound for $(\gamma_x,\dot\gamma_x)$. Indeed, the coercivity for \eqref{E-L}, as well as the fact that 
$x\mapsto \nabla R(t,x)$ grows linearly with a constant only depending on $t$, implies the existence of a locally bounded constant $K(t)$ such that 
\begin{equation}
\vert(\gamma_x(s),\dot\gamma_x(s))\vert\leq K(t)(1+\vert x\vert), \qquad \text{for all $s\in [0,t]$}.
\end{equation}
This implies,  modifying the constant $K(t)$ if necessary, thanks to \eqref{e2.2}:
$$
\vert u(t,x) \vert \leq K(t)(1+\vert x\vert^2).
$$
Moreover, from \eqref{e2.4}, we have
\beq
\label{e2.5}
\vert \nabla u(t,x)\vert \leq K(t)(1+\vert x\vert).
\eeq

\noindent{\bf Semi-convexity.}  Because $u$ is three times differentiable in $x$, then $\partial_tu$ is locally $W^{2,\infty}$ in $x$ and equation \eqref{e1.1} may be differentiated twice with respect to $x$. So, let $e$ be any unit vector, we have
$$
\begin{array}{rl}
\partial_t(\partial_{ee}u) &= 2 | \nabla \partial_e u|^2+2\nabla u.\nabla (\partial_{ee}u)+\partial_{ee}R(t,x)\\
&\geq
2 | \partial_{ee} u|^2+2\nabla u.\nabla (\partial_{ee}u)+\partial_{ee}R(t,x).
\end{array}
$$
%
Because of \eqref{e2.5}, the curve $t\mapsto\gamma_x(t)$ becomes a characteristic curve for the equation 
$$
\partial_tv=2 v^2 + 2\nabla u(t,x).\nabla v(t,x)+\partial_{ee}R(t,x).
$$
Moreover, along this characteristics, we find
$$
\f{d}{dt} v(t,\gamma_x(t))=2 v(t,\gamma_x(t))^2+\partial_{ee}R(t,\gamma_x(t)).
$$
We deduce, thanks to  \fer{asuD2} and \fer{asrD2-t}, that
$$
\partial_{ee}u \geq - \max( 2 \, \underline L_1,\sqrt{\underline K_1}).
$$

 \noindent{\bf Bounds on $D^3u$.} 
 As for the third derivative, we set $v(t,x)=D^3u(t,x)$, the equation for $v$ - that we may obtain using differential quotients  - is 
 $$
 \partial_tv - 2\nabla v.\nabla u = S(t,x,v):=6v.D^2u+D^3R,
 $$
 where $\nabla v$  denotes the column of tensors $(\partial_1v,...,\partial_dv)$ and $v.D^2u$ denotes the column of matrices $(\partial_1D^2u.D^2u,...,\partial_dD^2u.D^2u)$. This is a linear equation with (thanks to the bound on $D^2u$) bounded coefficients. Thus,  local boundedness of
 $\Vert v(t,.)\Vert_{L^\infty(\RR^d)}$ holds, and this  concludes the proof of Theorem \ref{cauchy}.

 \section{Uniqueness: reduction to a differential system and properties of  the function $I$}
 \label{sec:red}
 
 The idea is   to change the constrained problem (\ref{HJ}) by the following slightly nonstandard differential system:
\begin{equation}
\label{system}
\begin{cases}
R \left(\overline x(t),I(t) \right)=0,& \hbox{for $t\in \RR^+$},\\
\dot{\overline x}(t) = \left( -D^2u \big(t, \bar x(t)\big) \right)^{-1} \nabla R\big(\bar x(t), I(t) \big),& \hbox{for $t\in\RR^+$},\\
\partial_t u = | \nabla  u |^2 +R(x,I),& \hbox{in $\RR^+\times \R^d$},
\end{cases}
\end{equation}
with initial conditions
\begin{equation}
\label{first}
\begin{array}{c}
I(0)=I_0,\qquad u(0,\cdot)=u_0(\cdot),\qquad \overline x(0)=\overline x_0,\\
 \hbox{such that \; $\max_x u_0(x)=u_0(\overline x_0)=0$ \; and \; $R(\overline x_0,I_0)=0$.} 
\end{array}
\end{equation}
Note that (\ref{system}) is really a differential system because the assumptions on $R$ imply that $I(t)$ can implicitely be expressed in terms of $\bar x(t)$. And it is slightly nonstandard because $\bar x$ solves an ODE whose nonlinearity depends on $u$. The precise statement is the following
\begin{theorem}
\label{t3.1}
Solving the constrained problem \fer{HJ} is equivalent to solving the initial value   ODE-PDE problem \fer{system}-\fer{first}.
\end{theorem}
\noindent{\sc Proof.} Let $(u,I)$ be a solution of \fer{HJ} with initial datum $(u_0,I_0)$, the function $I$ being continuous, and $u$ a solution of the Hamilton-Jacobi equation in the sense of \fer{dyn-p}. Theorem \ref{cauchy} is applicable, and yields a solution $u(t,x)$ which has at least three locally bounded spatial derivatives, locally uniformaly in time. Moreover, the $D^2u$ is bounded uniformly in time and in $x$, and finally the function $u(t,.)$ is strictly concave. This allows a lot.
\begin{itemize}
 \item There is, at each time, a unique $\overline x(t)$ maximising $u(t,.)$  over $\RR^d$. Thus the trivial identity
\beq
\label{e3.01}
\nabla u(t,\overline x(t))=0
\eeq
can (use differential quotients)  be differentiated with respect to $t$, to yield that (i) $\overline x(t)$ is locally $W^{1,\infty}$ and (ii) the (a priori less trivial) identity
\beq
\label{e3.2}
\partial_t(\nabla u)(t,\overline x(t))+D^2u(t,\overline x(t)).\frac{d\overline x}{dt}(t)=0.
\eeq
\item The function $u(t,x)$ has enough regularity so that we may take the gradient of \fer{HJ} with respect to $x$, and evaluate the result at $x=\overline x(t)$. Because of \eqref{e3.01} we have $D^2u(t,\overline x(t)).\nabla u(t,\overline x(t))=0$ and, because of \eqref{e3.2}, we have
 \beq
\label{e3.3}
-D^2u(t,\overline x(t)).\frac{d\overline x}{dt}(t)=\nabla R(\overline x(t),I(t)).
\eeq
\item  The last item to take into account is the constraint
$$
u(t,\overline x(t))=0,
$$
which we may (still with the use of differential quotients) differentiate with respect to time, in order to yield
$$
\partial_tu(t,\overline x(t))+\nabla u(t,\overline x(t)).\frac{d\overline x}{dt}(t)=0,
$$
thus entailing 
$$
\partial_tu(t,\overline x(t))=0.
$$
This yields, from \fer{HJ}, 
\beq
\label{e3.5}
R(\overline x(t),I(t))=0.
\eeq
\end{itemize}
Gathering \eqref{e3.5}, \fer{HJ} and \fer{e3.3} shows that the constrained problem implies \fer{system}. \\

We also prove that  regularity plus \fer{system} easily implies \fer{HJ}.  Let $(u,I)$ solve \fer{system}. The first line of \fer{HJ} derives immediately. To prove the second line, note that, thanks to Theorem \ref{cauchy} and the third line of \fer{system} we deduce that $u$ is strictly concave. It has hence, for all $t\in \R^+$, a unique strict maximum point, in the variable $x$, that we denote $\overline y(t)$. Following similar arguments as above, we obtain that
$$
\begin{cases}
\dot{\overline y}(t) = (- D^2 u \left( t, \overline y (t) \right)^{-1} \, \nabla R( \overline y(t), I(t) ), & \text{for $t\in \R^+$,}
\\
\overline y (0) =\overline x_0.&
\end{cases}
$$
Comparing this with the second line of \fer{system} and \fer{first} we obtain that $\overline x(t) = \overline y(t)$, for all $t\in \R^+$. Finally, evaluating the third line of \fer{system} at $\overline x(t)$ and using the first line of \fer{system}, we obtain that $\p_t u(t,\overline x(t)) =0$. This equality together with $\nabla u(t,\overline x(t))=0$ and \fer{first} implies that
$$
\max_x u(t,x)=u(t,\overline x(t))=0.
$$
and the proof of Theorem \ref{t3.1} is complete. \hfill$\square$

\medskip
\noindent{\bf Remark.} {\it What we have done here is nothing else than the derivation of the equation for $\overline x$ carried out in \cite{DJMP}. In particular, the equilibrium property $R(\overline x,I)=0$ would hold in a more general setting than here. The new point here is that we establish, in a mathematically rigorous fashion,  the equivalence between the initial problem \fer{HJ} and 
 the coupled system \fer{system}, and this equivalence holds because of all the differentiations with respect to $x$ and $t$  that we are allowed to make. }

\section{The proof of Theorem \ref{thm:uniq}}
We fix $T>0$. To prove that \fer{HJ} has a unique solution $(u,I)$ in $[0,T]\times \R^d$, it is enough to prove that 
 there exists a unique solution to \fer{system}--\fer{first}.
 \\
 

We prove this using the Banach fixed point Theorem in a small interval and then iterate. 

\subsection{The mapping $\Phi$ and its domain}
First, we define 
$$
\begin{array}{c}
\Omega = \left\{ x \, | \, R(x,0)>0 \right\}, \qquad  \Omega_0= \left\{ x \,| \, R(x,I_0) \geq 0 \right\}, \quad \text{and} \\
\mathcal A=\left\{ x(\cdot) \in \mathrm{C} \Big( [0,\da]; B \big(\overline x_0,r_\da \big) \Big)   \, \Big| \, x(0)=\overline x_0 \right\},
\end{array}
$$
where $B(z,r)$ is the ball of radius $r$ centered at $z$, and $\da$ is a positive constant such that
\beq
\label{c-da}
 \da < \min( \mu, c(T)),  \quad \text{with }  \quad \mu= \f{\min (2\overline L_1,\sqrt{\overline K_1})}{C_M} \,     d(\Omega_0 \,,\, \Omega),
\eeq
where $d(A,B)$ is the distance between the  sets $A$ and $B$. The constant $r_\da$ is given by
$$
r_\da = \f{C_M \da}{\min (2\overline L_1,\sqrt{\overline K_1})},
$$
the constant $C_M$ is chosen such that
\beq
\label{CM}
| \nabla R (x,I ) |\leq C_M, \qquad \text{in $\Omega\times  [0,I_M]$},
\eeq
 and $c(T)$, a constant depending only on $T$,  will be chosen later. Note that, by the choice of $\da$ and since $\overline x_0 \in \Omega_0$, we obtain that
$$
 B\left( \overline x_0,r_\da \right) \subset \Omega.
$$

Our theorem will be proved by the introduction of a mapping 
$$
\Phi :\mathcal A\to \mathcal A,\qquad  \Phi(x)=y.
$$
We will prove the following theorem. 
\begin{theorem}
\label{t4.1}
The mapping $\Phi$ is a strict contraction from $\mathcal A$ into itself.
\end{theorem}

To define $\Phi$, we first need to introduce some other mappings. Let $x(\cdot) \in \mathcal A$. We define $\mathcal I : \mathcal A \to  \mathrm{C} \left( [0,\da] ; [0,I_M] \right)$ such that
$$
R(x(t),\mathcal I[x](t))=0.
$$
From \fer{asrmax}, \fer{asrDi} and since $x(t) \in \Omega$ for all $t\in [0, \da]$, it follows that $\mathcal I$ can be defined in a unique way. Next, we define the following domain 
$$
\begin{array}{rl}
\mathcal B = & \left\{ v\in L^\infty \big([0,\da]; W^{3,\infty}_{\mathrm loc}(\R^d) \big)\cap W^{1,\infty} \big([0,\da]; L^{\infty}_{\mathrm loc}(\R^d) \big)
\, |\,
\right. 
\\ 
&- \max( 2 \, \underline L_1,\sqrt{\underline K_1}) \leq D^2 v \leq - \min (2\overline L_1,\sqrt{\overline K_1}),\quad
 \left. \|D^3 v\|_{L^\infty([0,\da]\times \R^d)}  \leq L_4(T) \right\},
\end{array}
$$
and the following mapping
$$
\begin{cases}
V: \mathrm{C} \left( [0,\da] ; [0,I_M] \right) \to \mathcal B\\
V(I)=v,
\end{cases}
$$
where $v$ solves
\beq
\label{defV}
\begin{cases}
\p_t v = | \nabla v |^2 + R(x, I),& \text{in $[0,\da]\times \R^d$},\\
v(0,x) =u_0(x), &\text{in $\R^d$}.
\end{cases}
\eeq
It follows from Theorem \ref{cauchy} that the above mapping is well-defined. \\
Finally, we introduce a last mapping:
$$
\begin{cases}
F: \mathrm{C} \left( [0,\da] ; [0,I_M] \right) \times \mathcal A \times \mathcal B \to \mathcal A,\\
F(I,x,v)=y,
\end{cases}
$$
where $y\in \mathcal A$ solves
$$
\begin{cases}
\dot{y}(t) = \left( -D^2 v(t, x(t) \right)^{-1} \nabla R (x(t),I(t) ),& \text{in $[0,\da]$,}\\
y(0)=x_0.
\end{cases}
$$
To prove that $F$ is well-defined, we must verify that $y(t)$ remains in $B\left(x_0,r_\da \right)$.
 We note that, since $v\in \mathcal B$, we have
$$
0<\f{1}{ \max( 2 \, \underline L_1,\sqrt{\underline K_1}) } \leq \left( -D^2 v(t, x(t)) \right)^{-1} \leq \f{1}{ \min (2\overline L_1,\sqrt{\overline K_1})}.
$$
We deduce, thanks to \fer{CM}, that
$$
y(t) \in B\left(x_0, r_\da \right).
$$

\medskip

We are now ready to define mapping $\Phi$:
$$
\begin{cases}
\Phi : \mathcal A \to \mathcal A,\\
\Phi(x) = F\big( \mathcal I(x), x, V(\mathcal I(x)) \big).
\end{cases}
$$
It follows from the above arguments that $\Phi$ is well-defined. 
\subsection{Proof of Theorem \ref{thm:uniq}}
In this section we will explain why uniqueness to \eqref{system} holds. One main technical lemma  will be stated, its proof will be postponed to a special section.

Let us prove that $\Phi$ is a contraction for $c(T)$ (and hence $\da$) small enough. To this end, we first prove that $\mathcal I$ is Lipschitz:
\beq
\label{LipI}
| \mathcal I (x_1)- \mathcal I (x_2) | \leq C \| x_1 -x_2 \|_{L^\infty([0,\da])},\qquad \text{for all }x_1,\,x_2\, \in
 \mathrm{C} \Big( [0,\da] ; B \big( x_0, r_\da \big) \Big).
\eeq
We have indeed
$$
R(x_1,  \mathcal I (x_1)) = R(x_2,  \mathcal I (x_2))=0.
$$
It follows that
$$
R(x_1,  \mathcal I (x_1)) - R(x_2,I(x_1))= R(x_2,  \mathcal I (x_2)) - R(x_2,I(x_1)),
$$
and thus
$$
\nabla R \left( c x_1+(1-c) x_2,I(x_1) \right) \cdot (x_2-x_1) = \f{\p R}{\p I} (x_2,J) (I(x_2) - I(x_1)),
$$
with $c\in (0,1)$ and $J \in (I(x_1),I(x_2))$. Finally, using  \fer{asrDi}, \fer{CM} and the fact that $x_1(t),\, x_2(t) \in \Omega$ for all $t \in [0,\da]$ we obtain \fer{LipI} with $C=\f{C_M}{\overline K_2}$.

Next, we have that  $V:\mathrm{C} \left( [0,\da ]  ; [0,I_M]  \right) \to \mathcal B$ is also Lipschitz.
\begin{lemma}
\label{lem:LipV}
Let $I_1,\, I_2 \in \mathrm{C} \left( [0,\da] ; [0,I_M] \right)$. Then
\beq
\label{LipV}
\|V(I_1)-V(I_2)\|_{W^{2,\infty}([0,\da]\times \R^d)} \leq C \| I_1 -I_2 \|_{L^\infty([0,\da])} \da.
\eeq
\end{lemma}
This is a nontrivial lemma, whose proof will be given in the next section.
\medskip

From \fer{asrD2} and \fer{asr23}, and since $x_1, x_2 \in \mathcal A$ and $V_1, \, V_2 \in \mathcal B$, we deduce that $F: \mathrm{C} \left( [0,\da] ; [0,I_M] \right) \times \mathcal A \times \mathcal B \to \mathcal A$ is Lipschitz with respect to all the variables with Lipschitz constant $C\da$:
\beq
\label{LipF}
\begin{array}{rl}
\| F(I_1,x_1,V_1)-F(I_2,x_2,V_2) \|_{L^\infty([0, \da])}  &	\leq C \da
\Big[   \|x_1 - x_2�\|_ {L^\infty([0, \da])} + \| I_1-I_2 \|_{L^\infty([0, \da])}\\
&
+ \|V_1- V_2 \|_{W^{2,\infty}([0,\da]\times \R^d)} \Big].
\end{array}
\eeq
Finally, we conclude from \fer{LipI}, \fer{LipV} and \fer{LipF} that $\Phi : \mathcal A \to \mathcal A$ is a Lipschitz mapping with a Lipschitz constant 
$C\da$:
$$
\| \Phi(x_1) -  \Phi(x_2)  \|_{L^\infty([0, \da])} \leq C \da \| x_1 -x_2 \|_{L^\infty([0, \da])}.
$$
Choosing $c(T)$ (and hence $\da$) small enough, we deduce that $\Phi$ is a contraction.\\
 We deduce from the Banach fixed point Theorem that $\Phi$ has a unique fixed point and consequently \fer{system}--\fer{first} has a unique solution for $t\in [0,\da]$.
\\

 To prove that \fer{system}--\fer{first} has a unique solution in $[0,T]$ we iterate the above procedure $K=\lceil \f T \da \rceil$ times. Let $2\leq i\leq K$
  and $(x,I,u)$ be the unique solution of \fer{system}-\fer{first} for $t\in [0, (i-1)\da]$.
  Then, at the $i$-th step, we consider the same mapping $\Phi$ but as the initial data we choose
\beq
\label{ini:i}
 \widetilde x_0 =x((i-1)\da), \qquad \widetilde u_0 (\cdot)= u((i-1)\da, \cdot), \qquad \widetilde I_0 = I ((i-1)\da).
 \eeq
We claim that these initial conditions satisfy
\beq
\label{cond}
\begin{array}{c}\\
\widetilde x_0 \in \Omega_0, \qquad \widetilde u_0 \in   W^{3,\infty}_{\mathrm loc}(\R^d), \quad \text{ with 
$- \max( 2 \, \underline L_1,\sqrt{\underline K_1}) \leq D^2 \widetilde u_0 \leq - \min (2\overline L_1,\sqrt{\overline K_1})$}
\\  \|D^3 \widetilde u_0\|_{L^\infty( \R^d)}  \leq L_4((i-1)\da), \qquad
\max_x \widetilde u_0(x)= \widetilde x_0, \qquad R(\widetilde x_0, \widetilde I_0)=0.
\end{array}
\eeq
\begin{lemma}
\label{lem:comp}
Let $(x,I,u)$ be the unique solution of \fer{system}-\fer{first} in $[0,\tau]$. Then, $I(t)$ is increasing with respect to $t$ in $[0,\tau]$ and
$$
\begin{array}{c}
\overline x(t) \in  \Omega_0, \quad u(t,\cdot) \in   W^{3,\infty}_{\mathrm loc}(\R^d), \quad \text{ with \;
$- \max( 2 \, \underline L_1,\sqrt{\underline K_1}) \leq D^2 u(t,x) \leq - \min (2\overline L_1,\sqrt{\overline K_1})$,}\\
\|D^3 u\|_{L^\infty( [0,\tau]\times \R^d)}  \leq L_4(\tau), \quad
\max_x u(t,x)=\overline x(t), \quad R(\overline x(t), I(t))=0, , \quad \text{for all $t\in [0,\tau]$ and $x\in \R^d$}.
\end{array}
$$
\end{lemma}
\proof
The regularity estimates on $u$ are immediate from Section \ref{sec:pre}. The two last claims also follow immediately from \fer{system}-\fer{first} and the arguments in Section \ref{sec:red}. We only prove that   $I(t)$ is increasing with respect to $t$  and $\overline x(t)\in  \Omega_0$. 
 
 Differentiating the first line of \fer{system} with respect to $t$ we obtain
 $$
 \nabla R(\overline x(t), I(t)) \cdot \dot{\overline x}(t) + \f {\p }{\p I} R(\overline x(t), I(t)) \, \dot{I}(t)=0.
 $$
 Moreover, multiplying the second line of \fer{system} by $\nabla R$ we obtain
 $$
 \nabla R(\overline x(t), I(t)) \cdot \dot{\overline x}(t) =  \nabla R(\overline x(t), I(t)) \left(-D^2 u(t,\overline x(t) \right)^{-1}\nabla R(\overline x(t), I(t)) \geq 0.
 $$
 Combining the above lines we obtain
 $$
 \f {\p }{\p I} R(\overline x(t), I(t)) \, \dot{I}(t) \leq 0, \qquad \text{for $t\in [0,\tau]$},
 $$
 and hence, thanks to \fer{asrDi} we deduce
 $$
 \dot{I}(t) \geq 0, \qquad \text{for $t\in [0,\tau]$}.
 $$
 Therefore, $I(t)$ is increasing with respect to $t$ and in particular
 $$
 I(t) \geq I_0, \qquad \text{for $t\in [0,\tau]$}.
 $$
 Consequently, from the first line of \fer{system} and \fer{asrDi}, we obtain that
 $$
 R(\overline x(t), I_0) >0, \qquad \text{for $t\in [0,\tau]$}.
 $$ 
 It follows that
 $$
 \overline x(t) \in \Omega_0, \qquad \text{for $t\in [0,\tau]$}.
 $$
which concludes the proof. \qed
It is then immediate that the initial data given by \fer{ini:i} verify \fer{cond}. One can verify that the above conditions are the only properties that we have used to prove that $\Phi$ is well-defined and a contraction.  Therefore, one can apply again the Banach fixed point Theorem and deduce that there exists a unique solution of  \fer{system}-\fer{first} for $t\in [(i-1)\da, i\da]$.

\subsection{Proof of Lemma \ref{lem:LipV}}
(i) We first prove that 
\beq
\label{LipV1}
\|V(I_1)-V(I_2)\|_{L^{\infty}([0,\da]\times \R^d)} \leq C \| I_1 -I_2 \|_{L^\infty([0,\da])} \da.
\eeq
Let $v_1= V(I_1)$ and $v_2=V(I_2)$, and $r=v_2-v_1$. From \fer{defV} we obtain
\beq
\label{eq:r}
\begin{cases}
\p_t r= \left( \nabla v_1 + \nabla v_2 \right) \cdot \nabla r +R(x,I_2) -R(x,I_1),& \text{in $[0,\da]\times \R^d$}\\
r(0,x)=0,& \text{for all $x\in \R^d$.}
\end{cases}
\eeq
Note that the above equation has a unique classical solution which can be computed by the method of characteristics. The characteristics verify
\beq
\label{gam}
\dot{\gamma }(t) = - \nabla v_1 (t,\gamma)-\nabla v_2   (t,\gamma).
\eeq
For any $(t_1,x_1) \in [0,\da] \times \R^d$, there exists a unique characteristic curve $\gamma$ which verifies $\gamma (t_1)=x_1$. Moreover, this characteristic curve is defined in $[0,t_1]$.\\
The local existence and uniqueness of the characteristic curve $\gamma$ is derived from the Cauchy-Lipschitz Theorem and the fact that $\nabla v_i$ is Lipschitz with respect to $x$, for $i=1,2$. The latter property derives from the fact that $2 \underline L_1\leq D^2v_i \leq -2 \overline L_1$, for $i=1,2$.\\ 
To prove that the characteristic curve $\gamma$ is defined in $[0,t_1]$ we must prove that $\gamma$ remains bounded in this interval. However, this property follows using \fer{gam} and the fact that, since $v_1,\, v_2 \in \mathcal B$,
$$
|\nabla v_i(\gamma)| \leq C_1 |\gamma|+C_2, \qquad \text{for $i=1,2$}.
$$

We are now ready to prove \fer{LipV}. To this end, we multiply \fer{gam} by $\nabla r (t,\gamma)$ to obtain
$$
\nabla r (t,\gamma) \cdot \dot{\gamma }(t) = -( \nabla v_1 (t,\gamma)+\nabla v_2   (t,\gamma))\cdot \nabla r (t,\gamma).
$$
Combining this with \fer{eq:r} we deduce that
$$
\nabla r (t,\gamma(t))\cdot \dot{\gamma }(t) +\p_t r (t,\gamma(t)) = R(x,I_2) -R(x,I_1).
$$
We integrate this between $0$ and $t_1$ to find
$$
r (t_1,\gamma(t_1)) =\int_0^{t_1} R(x,I_2(\tau)) -R(x,I_1(\tau))d \tau.
$$
It follows, using \fer{asrDi}, that
$$
|r (t_1,\gamma(t_1)) | \leq  \underline K_2 |I_1-I_2| t_1,
$$
and hence \fer{LipV1}.\\

(ii) Next, we prove that
\beq
\label{bnr}
\|\nabla r\|_{L^{\infty}([0,\da]\times \R^d)} \leq C \| I_1 -I_2 \|_{L^\infty([0,\da])} \da.
\eeq
We differentiate \fer{eq:r} in the direction $e_i$, for $i=1,\cdots,d$, to obtain
$$
\p_t r_{(i)}= \left( \nabla v_{1,(i)} + \nabla v_{2,(i)} \right) \cdot \nabla r+\left( \nabla v_{1} + \nabla v_{2} \right) \cdot \nabla r_{(i)}  +R_{(i)}(x,I_2) -R_{(i)}(x,I_1),
$$
with the notation $f_{(i)}= \nabla  f\cdot e_i$. We multiply by $r_{(i)}$, sum over $i$ and divide by $|\nabla r|$ to obtain
$$
\p_t |\nabla r|  \leq \left( \nabla v_1 +\nabla v_2 \right) \cdot \nabla |\nabla r| + \sum_{i=1}^n \left| \nabla v_{1,(i)}+ \nabla v_{2,(i)} \right| |r_{(i)}| + |\nabla R (x,I_2) -\nabla R(x,I_1) |.
$$
Since $D^2 v_j$, for $j=1,2$, is bounded, and using \fer{asr23} we deduce that
 $$
\p_t |\nabla r|  \leq \left( \nabla v_1 +\nabla v_2 \right) \cdot \nabla |\nabla r| + C| \nabla r| + K_3 | I_2 - I_1 |.
$$
The characteristic curves corresponding to the above equation verify again \fer{gam}. We multiply \fer{gam} by $\nabla  |\nabla r|(t,\gamma(t))$ to obtain
$$
\nabla |\nabla r|  (t,\gamma) \cdot \dot{\gamma }(t) = -( \nabla v_1 (t,\gamma)+\nabla v_2   (t,\gamma))\cdot \nabla |\nabla r|  (t,\gamma).
$$
Combining the above equations we obtain
$$
\f{d}{dt} |\nabla r|  (t,\gamma(t)) \leq C| \nabla r|  (t,\gamma(t)) + K_3 \| I_2 - I_1 \|_{L^\infty([0,\da])}.
$$
It follows that
$$
|\nabla r (t,x)| \leq \left(e^{Ct}-1\right) \f{K_3}{C} \|I_2 - I_1\|_{L^\infty([0,\da])}.
$$ 
Hence \fer{bnr}, modifying the constant $C$ if necessary. 
\\

(iii) Finally we prove that
\beq 
\label{bJr}
\|D^2 r\|_{L^{\infty}([0,\da]\times \R^d)} \leq C \| I_1 -I_2 \|_{L^\infty([0,\da])} \da.
\eeq
Note that at every point $(t,x)\in \R^+\times \R^d$, we can write mixed derivatives of the form $r_{\xi\eta}$ in terms of pure derivatives:
\begin{equation}
\label{e3.1}
r_{\xi\eta}=\frac12(\partial_{\xi+\eta,\xi+\eta}^2r-\partial_{\xi\xi}^2r-\partial_{\eta\eta}^2r).
\end{equation}
This implies the existence of $\xi$ on the unit sphere of $\R^d$ such that
$$
\Vert D^2r\Vert_{L^\infty([0,\delta]\times\R^d)}\leq\frac32\Vert r_{\xi\xi}\Vert_{L^\infty([0,\delta]\times\R^d)}.
$$

We differentiate \fer{eq:r} twice in the direction of  $\xi$ and obtain
$$
\p_t r_{\xi\xi}= \left( \nabla v_{1,\xi\xi} + \nabla v_{2,\xi\xi} \right) \cdot \nabla r+
2 \left( \nabla v_{1,\xi} + \nabla v_{2,\xi} \right) \cdot \nabla r_{\xi}
+\left( \nabla v_{1} + \nabla v_{2} \right) \cdot \nabla r_{\xi\xi}  +R_{\xi\xi}(x,I_2) -R_{\xi\xi}(x,I_1).
$$
Using the above arguments, the fact that $D^2 v_i$ and $D^3 v_i$ are bounded, \fer{asr23} and \fer{bnr} we deduce that
$$
\p_t |r_{\xi\xi}| \leq C \| I_1 -I_2 \|_{L^\infty([0,\da])} \da +C\Vert r_{\xi\xi}\Vert_{L^\infty([0,\delta]\times\R^d)}+\left( \nabla v_{1} + \nabla v_{2} \right) \cdot \nabla |r_{\xi\xi}| + K_3 | I_2 - I_1 |.
$$
Next we use the characteristic curves as previously. We multiply \fer{gam} by $\nabla |r_{\xi \xi}(t,\gamma(t))|$ and obtain
$$
 \nabla |r_{\xi \xi}  (t,\gamma(t))| \cdot{\gamma }(t) = -( \nabla v_1 (t,\gamma)+\nabla v_2   (t,\gamma))\cdot \nabla |r_{\xi \xi}  (t,\gamma)|.
$$
We combine the above equations to obtain
$$
\f{d}{dt} |r_{\xi\xi}(t,\gamma(t))| \leq C \| I_1 -I_2 \|_{L^\infty([0,\da])} +C \Vert r_{\xi\xi}\Vert_{L^\infty([0,\delta]\times\R^d)}.
$$
We conclude that
$$
|r_{\xi \xi} (t,x)| \leq  Ct(  \|I_2 - I_1\|_{L^\infty([0,\da])}+\Vert r_{\xi\xi}\Vert_{L^\infty([0,\delta]\times\R^d)}).
$$ 
Restricting $c(T)$ (and hence $\delta$) if necessary we get
$$
\Vert r_{\xi\xi}\Vert_{L^\infty([0,\delta]\times\R^d)}\leq\frac{C\delta}{1-C\delta} \|I_2 - I_1\|_{L^\infty([0,\da])},
$$
hence \fer{bJr}. \qed

\section{An example with quadratic $u_0$ and $R$}

Let us study the (instructive) example  of a quadratic equation \eqref{HJ}. In other words we choose
$$
\left\{
\begin{array}{rl}
u_0(x)&=-\frac{1}{2}A_0x \cdot x,\\
R(x,I)&=-\frac12A_1x \cdot x+ b \cdot x+I_0-I,
\end{array}
\right.
$$
where $A_0$ and $A_1$ are positive definite matrices and $b\in\RR^d$. 
The Euler-Lagrange equation for the dynamic programming principle writes
\beq
\label{e5.10}
\left\{
\begin{array}{rll}
-\ddot\gamma+2A_1\gamma=&2b,\\
-\dot\gamma(0)+2A_0\gamma(0)=&0,\  \   \  \gamma(t)=x.
\end{array}
\right.
\eeq
In order to translate the equivalent system \fer{system}, we need $\nabla R(I,x)$ and $D^2u(t,x)$. The first quantity is easily obtained:
$$
\nabla R(I,x)=-A_1x+b.
$$
We then solve \fer{e5.10}, differentiate with respect to $x$ and denote by $\Gamma$ the differential of $\gamma$ - the solution of \fer{e5.10} with respect to $x$, to obtain (after an elementary but tedious computation)
\beq
\label{e5.12}
\begin{array}{rll}
\Gamma(s)=&e^{s\sqrt{2A_1}}B(s)B(t)^{-1}e^{-t\sqrt{2A_1}},\\
B(s)=&I_d+e^{-2s\sqrt{2A_1}}(\sqrt{2A_1}+2A_0)^{-1}(\sqrt{2A_1}-2A_0).
\end{array}
\eeq
Notice that the result was expected because $\gamma$ is linear function of $x$, and that $B(s)$ is invertible for $0\leq s\leq t$, and the norms of $B$ and $B^{-1}$ are bounded uniformly in $s$ and $t$. And so we have
\beq
\label{e5.14}
D^2u(t,x)=-A_0 \Gamma(0)^2-\int_0^t\biggl(\frac{\dot\Gamma(s)^2}2+A_1\Gamma(s).\Gamma(s)\biggl)ds:=-C(t),
\eeq
and the matrix $C(t)$ is bounded away from $0$ and $+\infty$, uniformly with respect to $t$. This, by the way, is not easily seen on the formula \fer{e5.14}; the proof of Theorem \ref{cauchy} is 
still what one should use here. 
The equation for $\overline x(t)$ and $I(t)$ is thus straightforward:
$$
\left\{
\begin{array}{rll}
\dot{\overline x}(t)=&C(t)^{-1}(-A_1\overline x(t)+b), \\
I(t)=&I_0-\di\frac12A_1\overline x(t) \cdot \overline x(t)+b\cdot\overline x(t).
\end{array}
\right.
$$
Thus, existence and uniqueness of $(\overline x,I)$ is straightforward.
For large $t$, examination of \fer{e5.14} and \fer{e5.12}
yields 
$$
D^2u(t,x)\sim_{t\to+\infty}-\sqrt{\frac{A_1}{2}},
$$
and 
$$
\left\{
\begin{array}{rll}
\di\lim_{t\to+\infty}{\overline x}(t)=&A_1^{-1}b,\\
\di\lim_{t\to+\infty}I(t)=&I_0+\di\frac12A_1^{-1}b\cdot b.
\end{array}
\right.
$$
This is consistent with the known behaviour of $I(t)$ and $\overline x(t)$, as well as a closed form of the competition increase $I(t)-I_0$.

%
%

\section*{Acknowledgements}
The authors thank G.  Barles and B. Perthame for interesting comments on an earlier version of this paper, which had a great influence on the present form of this work. They also thank M. Bardi for pointing out the possible tractability of the quadratic case.
S. Mirrahimi was partially funded by the ANR projects KIBORD ANR-13-BS01-0004 and MODEVOL ANR-13-JS01-0009. 
J.-M. Roquejoffre  was supported by  the European Union's Seventh Framework Programme (FP/2007-2013) / ERC Grant
Agreement n. 321186 - ReaDi - ``Reaction-Diffusion Equations, Propagation and Modelling'' held by Henri Berestycki,
as well as the ANR project  NONLOCAL ANR-14-CE25-0013.  Both authors thank  the Labex CIMI for  a PDE-probability winter quarter in Toulouse, which provided a stimulating scientific environment to this project.

 \footnotesize{

}

\begin{thebibliography}{100}


\bibitem{Barles} {\sc G. Barles,} {\it Solutions de viscosit\'e des équations de Hamilton-Jacobi}, Math\'ematiques et Applications, Springer, 1994.

\bibitem{GB.LE.PS:90}
{\sc G.~Barles, L.~C. Evans, and P.~E. Souganidis},
{\it Wavefront propagation for reaction-diffusion systems of {PDE}},
 Duke Math. J., {\bf 61(3)} (1990) pp. 835--858.

\bibitem{BMP}
{\sc G.~Barles, S.~Mirrahimi,  B.~Perthame},
{\it Concentration in {L}otka-{V}olterra parabolic or integral equations:
  a general convergence result.}
 Methods Appl. Anal., {\bf 16(3)} (2009) pp.321--340.
 
 
 \bibitem{NC.phd}
{\sc N.~Champagnat},
{\it  Mathematical study of stochastic models of evolution belonging to the ecological theory of adaptive dynamics}, 
Ph.D. Thesis, University of Nanterre (Paris 10), 2004.
 

\bibitem{NC.RF.SM:06}
{\sc N.~Champagnat, R.~Ferri\`ere,  S.~M\'el\'eard,}
{\it  Unifying evolutionary dynamics: From individual stochastic processes
  to macroscopic models}, 
 Theoretical Population Biology, {\bf 69(3)} (2006) pp. 297--321.
 
\bibitem{CFM}
{\sc N.~Champagnat, R.~Ferri\`ere,  S.~M\'el\'eard},
 {\it Individual-based probabilistic models of adaptive evolution and
  various scaling approximations}, {\bf 59} (2008) Progress in Probability, Birkha\"{u}ser.

\bibitem{CJ}
{\sc N.~Champagnat, P.-E. Jabin},
{\it The evolutionary limit for models of populations interacting competitively via several resources,}
Journal of Differential Equations, {\bf 261} (2011) pp.179--195.
%


\bibitem{OD:04}
{\sc O.~Diekmann},
{\it Beginner's guide to adaptive dynamics,}
 Banach Center Publications {\bf 63} (2004) pp. 47--86.


\bibitem{DJMP}
{\sc O.~Diekmann, P.-E. Jabin, S.~Mischler,  B.~Perthame},
{\it The dynamics of adaptation: an illuminating example and a
  {H}amilton-{J}acobi approach},
 Th. Pop. Biol., {\bf 67(4)} (2005) pp. 257--271.
 
 
 \bibitem{LE.PS:89}
{\sc L.~C. Evans and P.~E. Souganidis},
{\it  A {PDE} approach to geometric optics for certain semilinear parabolic
  equations},
 Indiana Univ. Math. J., {\bf 38(1)} (1989) pp. 141--172.

\bibitem{MFb:85}
{\sc M.~Freidlin},
{\it Functional integration and partial differential equations},
Annals of Mathematics Studies,    {\bf   109} (1985) Princeton University Press, Princeton, NJ.

\bibitem{MF:85}
{\sc M.~Freidlin},
{\it Limit theorems for large deviations and reaction-diffusion equations},
 The Annals of Probability, {\bf 13(3)} (1985) pp. 639--675.
 
 
 \bibitem{FIL1} {\sc Y. Fujita,  H. Ishii, P. Loreti}, {\it  Asymptotic solutions of viscous Hamilton-Jacobi equations with Ornstein-
Uhlenbeck operator,} Comm. Partial Differential Equations, {\bf 31} (2006) pp. 827-848.
\bibitem{FIL2} {\sc Y. Fujita, H. Ishii, P.Loreti,} {\it Asymptotic solutions of Hamilton-Jacobi equations in Euclidean
$n$ space},
Indiana Univ. Math. J., {\bf 55} (2006) pp. 1671-1700.

\bibitem{SG.EK:98}
{\sc S.~A.~H. Geritz, E.~Kisdi, G.~M\'eszena, J.~A.~J. Metz.},
{\it Evolutionary singular strategies and the adaptive growth and
  branching of the evolutionary tree},
 Evolutionary Ecology, {\bf 12} (1998) pp. 35--57.

\bibitem{LMP} {\sc A. Lorz, S. Mirrahimi, B Perthame}, {\it Dirac mass dynamics in a multidimensional nonlocal parabolic equation,} Communications in Partial Differential Equations, {\bf 36} (2011) pp.1071--1098.



\bibitem{Mthese} {\sc S. Mirrahimi}, {\it Ph\'enom\`enes de concentration dans certaines EDPs issues de la biologie}, PhD thesis, Universit\'e Pierre et Marie Curie, 2011,  http://www.math.univ-toulouse.fr/~smirrahi/manuscrit.pdf.

\bibitem{SM.BP:15}
{\sc S.~Mirrahimi,  B.~Perthame},
{\it Asymptotic analysis of a selection model with space},
 To appear in  J. Math. Pures App..




\bibitem{MR0}  {\sc S. Mirrahimi, J.-M. Roquejoffre}, {\it Uniqueness in a class of Hamilton-Jacobi equations with constraints,} C. R. Acad. Sci. Paris, Ser. I (2015), http://dx.doi.org/10.1016/j.crma.2015.03.005.


\bibitem{MR2}  {\sc S. Mirrahimi, J.-M. Roquejoffre}, {\it Approximation of solutions of selection-mutation models and error estimates,} in preparation.
\bibitem{BP} {\sc  B. Perthame, G. Barles}, {\it Dirac concentrations in Lotka-Volterra parabolic PDEs,} Indiana Univ. Math. J. {\bf 57} (2008) pp. 3275--3301.


\bibitem{Raoulphd}
{\sc G.~Raoul},
{\it Etude qualitative et num\'{e}rique d'\'{e}quations aux
  d\'{e}riv\'{e}es partielles issues des sciences de la nature},
PhD thesis, ENS Cachan, 2009.

\end{thebibliography}
\end{document}